\documentclass[11pt]{amsart}
\usepackage{amsfonts,amssymb,amsmath,a4wide,xypic,remreset}

\newcommand{\Zm}{\ensuremath{\mathbb{Z}}}

\newcommand{\Rm}{\ensuremath{\mathbb{R}}}

\newcommand{\Tm}{\ensuremath{\mathbb{T}}}

\newcommand{\Nm}{\ensuremath{\mathbb{N}}}

\def\lto{\longrightarrow}
\def\lmto{\longmapsto}

\def\leq{\leqslant}
\def\geq{\geqslant}

\newcommand{\mF}{\ensuremath{\mathcal{F}}}

\author{Patrick Bernard}
\title{
Hamiltonian Systems: Stability and Instability Theory}

\begin{document}

\maketitle

The solar system has long appeared to astronomers and mathematicians
as a model of stability.
On the other hand, statistical mechanics relies  on the assumption
that large assemblies of particles form highly unstable systems
(at the microscopic scale).
Yet all these physical situations are described, at least to a certain
degree of approximation, by Hamiltonian systems.

One may hope that  Hamiltonian systems can be classified 
in two different categories, stable and unstable ones.
However,
the situation  is much more complicated and  both stable and 
unstable behaviors cohabit in typical systems.
Even our examples are not perfect  paradigms of stability and instability. 
Indeed, it is now clear from numerical as well as theoretical points
of views that some instability is present over long time-scales
in the solar systems, so that for example future collisions 
between planets can not be completely ruled out in view
of our present understanding.
On the other hand, unexpected patterns of stability have been 
discovered in systems involving a large number of particles.

Understanding the impact of stable and unstable effects 
in Hamiltonian systems 
has been considered since Poincar\'e as 
 one of the most important questions
in dynamical systems. 
In the present text, we will discuss model Hamiltonian systems
of the form
$$
H_{\epsilon}(q,p)=h(p)+\epsilon G_{\epsilon}(q,p)
$$
where $(q,p)\in \Tm^d\times U$,
with $U$ a bounded open subset of $\Rm^n$.
Recall that the equations of motion are 
\begin{equation}\label{eqnq}
\dot q(t)=\partial_ph(p)+\epsilon \partial_pG_{\epsilon}(q,p)
\end{equation}
\begin{equation}\label{eqnp}
\dot p(t)=-\epsilon \partial_qG_{\epsilon}(q,p).
\end{equation}
The textbook \cite{Ar} is a good general introduction
on Hamiltonian systems.
We will always denote by $\omega(p)$ the frequency map 
$\partial_ph(p)$, which plays a crucial role.
Here, as is obvious in (\ref{eqnp}), the action variables $p$
are preserved under the evolution in the unperturbed case
$\epsilon=0$.
We will try to explain what is known on the evolution
of these action variables for the perturbed system.
As we will see,
in many situations, these variables are extremely stable.
For example, KAM theorem implies that, for a positive measure
of initial conditions  $(q_0,p_0)$ the trajectory
$(q(t),p(t))$ satisfies
$\|p(t)-p(0)\|\leq C\epsilon$ for all times.
Examples show that some initial conditions may lead to unstable
trajectories, that is trajectories such that 
$\|p(t)-p(0)\|\geq 1/C$ for some $t$
(depending on $\epsilon$)  and some fixed constant $C$
independent of $\epsilon$.
However, this is, as we will see, possible only for very large
time $t$ (meaning that $t$ as a function of $\epsilon$
has to go to infinity very quickly when $\epsilon\lto 0$). 
The main questions here are to understand in what situation
instability is or is not possible, and what kind of evolutions
can have the actions variable $p$. 
Another important question is to estimate the speed
(as a function of the parameter $\epsilon$) of the evolutions
of $p$. 

\subsection{A convention}
We assume, unless otherwise stated, that the Hamiltonians
are real analytic. The norm $|H|$ of the Hamiltonian
$H$ is the uniform norm of its holomorphic extension to
a certain complex strip. We do not specify the width of
this strip.
Whenever we consider a family $H_{\epsilon}$, $F_{\epsilon}$\ldots
of Hamiltonians, we mean that the norm $|H_{\epsilon}|$
is bounded when $\epsilon \lto 0$.

\section{Averaging and exponential stability.}

The first  observation concerning the action variable
is that they should evolve at a speed of the order of 
$\epsilon$.
 However, averaging effects occur.
More precisely, in the equation
$\dot p(t)=-\epsilon \partial_q H_{\epsilon}(q(t),p(t))$,
the variable $q(t)$ is moving fast compared to $p(t)$.
If the evolution of $q(t)$ nicely fills the torus $\Tm^n$, 
it is tempting to think that the 
averaged equation
$$
\dot {\bar p}(t)=-\epsilon  \bar V_{\epsilon}(\bar p(t))
$$
should approximate accurately the actual behavior
of $p(t)$, where 
$$
\bar V_{\epsilon}(p):=\int_{\Tm^d}\partial _qH(q,p)dq.
$$
We have $\bar V\equiv 0$, which leads to think that the evolution
 should consist mainly of oscillations
of small amplitude with no large evolution.
This reasoning is limited by the presence of resonances.

\subsection{Frequencies}
A frequency 
 $\omega\in \Rm^d$ is said resonant if there exists 
$k\in \Zm^d_*  (=\Zm^d-\{0\})$ such that  $\langle k, \omega\rangle=0$.
The resonance module of 
 $\omega$,
$$
Z(\omega)=\{k\in \Zm^d / \langle k, \omega\rangle=0\}
$$
is a subgroup $\Zm^d$, we denote by
$R(\omega)$
the vector space generated by $Z(\omega)$ in  $\Rm^d$.
The order of resonance $r(\omega)$ is the dimension of 
$R(\omega)$.
The main examples of a resonances of order $r$ are the
frequencies $\omega=(\omega_1,0)$ where  $\omega_1\in \Rm^{d-r}$
is  non resonant.
This example is  universal.
Indeed,  if $\omega$ is a resonant frequency,
then there exists a matrix  $A\in Gl_d(\Zm)$
such that  $A\omega=(\omega_1,0)$, where  $\omega_1\in \Rm^{d-r}$ 
is not resonant.
The matrix $A$ can be seen as a diffeomorphism
of $\Tm^d$, which transports the constant vector-field $\omega$
to the constant vector-field $A\omega=(\omega_1,0)$.
It is useful to distinguish, among non resonant frequencies,
some which are sufficiently non-resonant.
A frequency $\omega\in \Rm^d$
is called Diophantine if there exists real constants $\gamma>0$ and 
 $\tau\geq d$
such that 
$$
|\langle \omega,k\rangle |
\geq \gamma \|k\|^{1-\tau}
$$
for each $k\in \Zm^d_*$.
Finally, a frequency is called resonant-Diophantine
if there exists a matrix $A\in  Gl_n(\Zm)$ such that
$\omega=A(\omega_1,0)$, where $\omega_1\in \Rm^{d_1}$
is a Diophantine frequency.

\subsection{Symplectic diffeomorphisms and Normal Forms}
An efficient mathematical method to take averaging effects
into account is the use of normal forms.
Normal form theory consists in finding new coordinates in which
the fast angles have been eliminated from the equations up 
to a small remainder. This is done exploiting the existence of a large
group of diffeomorphisms preserving the Hamiltonian structure of 
equations, called symplectic diffeomorphisms or 
canonical transformations.
We refer the reader to standard textbooks for these notions, 
for example
to \cite{Ar}.
An important point is that a symplectic diffeomorphism $\phi$
sends the trajectories of the Hamiltonian $H\circ \phi$
to the trajectories of the Hamiltonian $H$.
A Hamiltonian $N(q,p)$ is said in $R$-normal form,
where $R$ is a linear subspace of $\Rm^n$, if
$\partial_qN\in R$ for each $(q,p)$.
Let us give an illustrative result, taken from \cite{LMS}.
Note that this result is not sufficient to obtain uniform stability 
estimates, as in Nekhoroshev Theorem below. 
More precise normal form results are given in \cite{Ne} and \cite{Po}.

\subsection{Normal form theorem}
\begin{itshape}
Let $\omega_0=\omega(p_0)$ be a given Diophantine or resonant-Diophantine
frequency. Let us denote $B_r(p_0)$ the open ball of radius $r$ in $\Rm^d$
centered at $p_0$.
There exists a constant $a$ which depends only on $\omega$,
 and  constants $\epsilon_0>0$ and  $C>0$ 
such that the following holds:
For each $\epsilon<\epsilon_0$,
there exists an analytic  symplectic embedding
$\phi_{\epsilon}: \Tm^d\times B_{r(\epsilon)}\lto\Tm^d \times U $, 
which is
$\epsilon$-close to identity and such that  
$$
H_{\epsilon}\circ \phi_{\epsilon}(q,p) =
h(p)+\epsilon N_{\epsilon}(q,p)+\mu(\epsilon)F_{\epsilon}(q,p),
$$
where $N$ is in $R(\omega_0)$-normal form, $r(\epsilon)\geq \sqrt{\epsilon}$,
and 
 $\mu(\epsilon)\leq e^{-C\epsilon^{-a}}$.
\end{itshape}

This means that the motions with resonant initial conditions
are confined, up to small oscillations, in  the associated affine plane 
$p(0)+R(\omega(p(0))$ until they live the domain of the normal form,
or until time $\mu^{-1}(\epsilon)$.

\subsection{Geometry of resonances}\label{rescaled}
In view of the Normal Form Theorem, we are led to consider
the curves $P(\theta):\Rm\lto \Rm^d$ which satisfy
$$
P(\theta')-P(\theta)\in R(\omega(P(\theta)))
$$
for each $\theta$ and $\theta'$.
Indeed, it appears that these curves are the ones 
the action variables can follow on time-scales not 
involving the remainders of the normal forms.
Note that here the parameter $\theta$ is not the physical time.
Assuming that $P(\theta)$ is such a curve, 
we can define the affine space  
$$R:= P(0)+\cap_{\theta\in \Rm}R(\omega(P(\theta))).$$
We then  have $P(\theta)\in R$
for each $\theta$.
In addition, each point $P(\theta),\theta\in \Rm$
is a critical point of the restriction $h_{|R}$
of the unperturbed Hamiltonian $h$ to the affine space $R$.
It follows that the curve $P(\theta)$ has to be constant if
the unperturbed Hamiltonian satisfies the following hypothesis.

\subsection{Nekhoroshev steepness}
\begin{itshape}
We say that the unperturbed Hamiltonian $h$ is steep
if, for each affine subspace $\Lambda$ in $\Rm^d$,
the restriction $h_{|\Lambda}$
has only isolated critical points.
\end{itshape}

This formulation, due to  L. Niederman, \cite{Ni},  is much simpler
than the equivalent one first given by Nekhoroshev in \cite{Ne}.
It turns out that this condition, which was made natural by our heuristic
explanation, implies  stability over exponential time-scales
for all initial conditions, see \cite{Ne}.
We first need another condition.

\subsection{Kolmogorov non-degeneracy}\label{KND}
\begin{itshape}
We say that the unperturbed Hamiltonian $h$ is
non-degenerate in the sense of Kolmogorov
if it has non-degenerate Hessian at each point,
or equivalently if the frequency map $p\lmto \omega(p)$
is an immersion.
\end{itshape}

\subsection{Nekhoroshev stability theorem}\begin{itshape}
Assume that the unperturbed Hamiltonian does not have critical points 
($\omega(p)$ does not vanish), satisfies Nekhoroshev steepness
and Kolmogorov non-degeneracy conditions.
Then there exists constants $a>0$ and $b>0$, which depend only on $h$,
and constants $\epsilon_0>0$ and $C>0$
such that the following holds:
For  $\epsilon<\epsilon_0$, each trajectory
$(q_{\epsilon}(t),p_{\epsilon}(t))$ satisfies the estimate
$$
\|p_{\epsilon}(t)-p_{\epsilon}(0)\| \leq C{\epsilon}^{b}
$$
for all $t$ such that  $|t|\leq e^{C{\epsilon}^{-a}}$.
\end{itshape}

\subsection{Herman's example}\label{Herman}
In oder to illustrate the necessity of the condition of steepness,
let us consider the Hamiltonian 
$$
H_{\epsilon}(q_1,q_2,p_1,p_2)=p_1p_2+\epsilon V(q_1).
$$
with   $V:\Tm\lto \Rm$.
The associated equations are
$$
\dot p_2=0, \;\; \;\;\dot p_1=-V', \;\; \;\; \dot q_1=p_2,
  \;\; \;\;\dot q_2=p_1.
$$
The trajectories whose initial condition are subjected to
 $p_2(0)=0$  and  $V'(q_1(0))\neq 0$
satisfy
$$
p_1(t)=p_1(0)-t\epsilon V'(q_1(0)), 
\;\; \;\;p_2(t)=0,  \;\; \;\;q_1(t)= q_1(0).
$$
We see an evolution at speed $\epsilon$ of the action variable
$p_1$ contradicting the conclusion of Nekhoroshev theorem.
In this example, we have $R(\omega(p(t)))=\Rm\times\{0\}$,
and $h_{|\Rm\times\{0\}}\equiv 0$, so that the curve
$$
P(\theta)=(\theta,0)
$$
is indeed a curve of critical points of 
 $h_{|\Rm\times\{0\}}$.

\subsection{Genericity of steepness}
The condition of steepness is frequently satisfied.
In order to be more precise, we mention that, for 
$N\in \Nm$ large enough (how large depends on the dimension $d$),
steepness is a generic condition in the finite dimensional
space of polynomials of degree less that $N$.
Note in contrast that a quadratic Hamiltonian is steep
if and only if it is positive definite.
Finally, it is important to mention that convex Hamiltonians
$h$ with positive definite Hessian  are steep.
More generally, quasi-convex Hamiltonians are steep.
A function $h:U\lto \Rm$ is said quasi-convex
if, at each point, the restriction of its Hessian
to the kernel of its differential is positive definite. 

\subsection{The quasi-convex case}
It is interesting to be more precise about the values of $a$ and $b$
in Nekhoroshev Theorem.
We shall do so in the quasi-convex case, 
which is the most stable case, and where much more is known.
If $h$ is quasi-convex, one can take 
$$
a=b=\frac{1}{2d},
$$
as was proved by P. Lochak, see \cite{Lo}.
It is a question of active present research whether these exponents 
are optimal.
It now appears that this is almost so, and that the optimal
exponent $a$ should not be larger than $1/2(d-3)$.
That this exponent deteriorates as the dimension increases is
of course very natural in the perspective of statistical mechanics.
As a matter of fact, not only the exponent $a$, but also
the threshold $\epsilon_0$ of validity of Nekhoroshev Theorem
deteriorates with the dimension, as was noticed in \cite{BK}.

Another important fact was proved in \cite{Lo}:
in these expressions, the important value of $d$ is not the total
number of degrees of freedom, but the number of active degrees
of freedom. More precisely, resonant initial
condition are more stable than generic ones.
If $r$ is the order of resonance of a given initial condition,
then the number $d-r$ of fast angles can be substituted to
the total number of degrees of freedom for the computation
of the stability exponent. This phenomenon may
account for the surprising stability obtained numerically
by Fermi, Pasta and Ulam.

%\subsection{Normal form theorem}
%\begin{itshape}
%Let $\omega_0=\omega(p_0)$ be a given frequency, and let $r$ be the order 
%of resonance of $\omega_0$.
%Let $B_a(p_0)$ be the open ball of radius $a$ in $\Rm^d$
%centered at $p_0$.
%Then, there exists a constant $A>0$ and, for $\epsilon<1/A$,
%a symplectic embedding
%$\phi_{\epsilon}:\Tm^d\times B_{a(\epsilon)}(p_0)\lto \Tm^d\times \Rm^d$ 
%such that
%$$
%H_{\epsilon}\circ \phi_{\epsilon}=
%h(p)+N_{\epsilon}(q,p)+\mu(\epsilon)F_{\epsilon}(q,p),
%$$
%where 
%\begin{itemize}
%\item $N_{\epsilon}$ is a Hamiltonain in normal form, which means that 
%$\partial_qN\in R(\omega_0)$,
%\item  $|N_{\epsilon}|$ and  $|F_{\epsilon}|$ are bounded independantly 
%of $\epsilon$.
%\item
%$\epsilon^{2(d-r)}/A\leq a(\epsilon)\leq A \epsilon^{2(d-r)}$ 
%\item
%$\mu(\epsilon)\leq \exp\big( -\epsilon^{-2(d-r)}/A\big)$
%\end{itemize}
%\end{itshape}
%Unfortunately, 
% the constant $A$ can not be chosen uniformly with 
%respect to $\omega_0$. As a consequence, although quite illustrative,
%this statement is not sufficiently precise from a technical viewpoint
%in order to obtain uniform stability estimates.

\section{Permanent stability.} 
Many initial conditions
satisfy more than exponential stability: they are permanently stable.

\subsection{Kolmogorov Theorem}
\begin{itshape}
Assume that $h$ satisfies Kolmogorov non-degeneracy condition 
\ref{KND}. Then for each open subset $V\subset \Rm^d$ 
such that $\bar V\subset U$,
there exists $\epsilon_0>0$ such that,
for each $\epsilon<\epsilon_0$ there exists
\begin{itemize}
\item
a smooth symplectic embedding   
$\phi_{\epsilon}:\Tm^d\times V\lto
\Tm^d\times U
$, which is $\epsilon$-close to the identity,
\item 
a compact subset $F_{\epsilon}$ of $V$, whose relative measure
in $V$ is converging to $1$ as $\epsilon\lto 0$,
\end{itemize} 
such that the Hamiltonian system $H_{\epsilon}\circ \phi_{\epsilon}$
preserves the torus $\Tm^d\times \{p\}$ for each $p\in F_{\epsilon}$.
 \end{itshape}

The union 
$$
\mF_{\epsilon}=\phi_{\epsilon}(\Tm^d\times F_{\epsilon})
$$
of all the invariant tori has positive measure.
Its complement is usually an open dense subset of $\Tm^d\times U$.
All the orbits starting in this invariant set obviously undergo oscillations
of amplitude of the order of $\epsilon$ for all times.
It is worth mentioning that some energy surfaces may not
intersect the invariant set $\mF_{\epsilon}$.
This is illustrated in example \ref{Herman}, where the
surface of zero energy does not contain invariant tori.
The following condition guaranties the existence 
of invariant tori on each energy surface.

\subsection{Arnold non-degeneracy}
\begin{itshape}
The Hamiltonian $h$ is said to be non-degenerate 
in the sense of Arnold if it does not have critical points
and if the map 
$$
p\lmto \frac{\omega(p)}{\|\omega(p)\|}
$$
is a local diffeomorphism between each level set of $h$ and $S^{d-1}$.
This is equivalent to say that the function  
$(\lambda,p)\in \Rm\times U\lmto \lambda h(p)$
 has non-degenerate Hessian at 
each point of the form $(1,p)$.
\end{itshape}

\subsection{Arnold Theorem}
\begin{itshape}
If $h$ satisfies Arnold non-degeneracy condition, then the relative measure
of the set 
$\mF_{\epsilon}$ of invariant tori is converging to $1$ in each energy surface.
\end{itshape}

This theorem prevents ergodicity of the perturbed systems 
for the canonical invariant measure on its energy surface.
This may be  considered as a very disappointing result
for statistical mechanics, whose mathematical
foundation has often be  considered to be the Boltzmann hypothesis 
of ergodicity.
However, statistical mechanics is 
first of all a question of letting $d$ go to infinity,
and ergodicity  might not be such a crucial hypothesis,
see \cite{K}.

When $d=2$, The theorem of Arnold has particularly
strong consequences. Indeed, in this case,
the invariant tori cut the energy surfaces in small
connected components. The motion is then confined in these
connected components. As a consequence, we obtain permanent
stability for all initial conditions.

In higher dimension however, the complement of $\mF_{\epsilon}$
in each energy shell is usually a dense, connected open set.
There may exist orbits wandering in this large connected set,
although the speed of evolution of these orbits is limited
by Nekhoroshev theory. Understanding the dynamics 
in this open set is a very important and hard question.
It is the subject of the next section.

\subsection{Relaxed assumption}
For many applications, such as celestial mechanics,
the non-degeneracy conditions of Arnold or Kolmogorov
are not satisfied, or hard to check.
However, the existence of invariant tori has been proved under 
much milder  assumptions.
As a rule, Invariant tori exist in the perturbed systems
if the frequency map $p\lmto \omega(p)$ stably contains 
Diophantine vectors in its image.

\section{The Mechanism of Arnold}
Understanding instability is the subject of 
intense present research.
General methods of construction of interesting orbits,
as well as clever classes of examples are being developed.
These methods are exploring the limits of stability theory.
Here we shall only describe the foundational ideas of 
Arnold, see \cite{Ar}, where most of  the present activity 
finds its roots.
Although these ideas have some ambition of universality, they are best 
presented, like in \cite{Ar}, on an example.
We consider the quasi-convex Hamiltonian 
$$
H(q_1,q_2,q_3,p_1,p_2,p_3)=(p_1^2+p_2^2)/2-p_3+\epsilon \cos 2\pi q_2 +
\mu (\cos 2\pi q_2)(\cos2\pi  q_1+\cos 2\pi q_3).
$$
As we have seen, this system is typical of
the kind of Hamiltonians one gets after reduction
to resonant normal form. 
However, it is illuminating to consider 
 $\mu$
not as a function of $\epsilon$ but as an independent parameter.
This is an idea of Poincar\'e then followed by Arnold.
We shall expose the main steps of the proof of the following 
result.

\subsection{Theorem}
\begin{itshape}
Let us fix numbers 
$0<A<B$.
For each $\epsilon>0$, there exists a number 
 $\mu_0(\epsilon)$
such that, when  $0<\mu<\mu_0(\epsilon)$, 
there exists a trajectory
$$(q_1(t),q_2(t),p_1(t),p_2(t))$$
and a time $T>0$ (which depends on $\epsilon$ and $\mu$) such that 
$$p_1(0)\leq A\text{ , } p_1(T)\geq B.$$
\end{itshape}

\subsection{The truncated system}
Let us begin with some remarks about the 
truncated Hamiltonian obtained when $\mu=0$:
$$
H_0(q,p)=H_1(q_1,q_3,p_1,p_3)+H_2(q_2,p_2)=p_1^2/2 -p_3+
p_2^2/2+\epsilon\cos 2\pi q_2.
$$
This system is the uncoupled product of $H_1$ and 
of the pendulum described by $H_2$.
The variable $p_1$ is constant along motion, hence the Theorem
can not hold for $\mu=0$.

Recall that the point
$q_2=0, p_2=0$
is a hyperbolic fixed point of the pendulum
$H_2(q_2,p_2)=p_2^2/2+\epsilon \cos 2\pi q_2.
$
The stable and unstable manifolds of this integrable system
coincide, they from the energy level
$H_2=\epsilon$.
As a consequence, in the product system of Hamiltonian
$H_0=H_1+H_2$, there exists, in the zero energy level,  a one parameter
family  $T_{\omega}$ of invariant tori of dimension 2
$$
T_\omega=\{p_1=\omega, p_3= \omega^2/2+\epsilon,
q_2=0, p_2=0\}\subset \Tm^3\times\Rm^3.
$$
Each of these tori is hyperbolic
in the sense that it has a stable manifold of dimension 3 and an unstable
manifold of dimension 3, which are nothing but the liftings
of the stable and unstable manifolds of the hyperbolic fixed point of 
$H_2$.
Notice that these manifolds do not intersect transversally along
$T_{\omega}$.

When $\mu\neq 0$, the perturbation is chosen in such a way that
the tori $T_{\omega}$ are left invariant by the Hamiltonian flow.
\subsection{Splitting}
\begin{itshape}
For  $0<\mu<\mu_0(\epsilon)$,
the invariant tori $T_{\omega}$ still have stable and unstable manifolds
of dimension 3. 
These stable and unstable manifolds intersect 
transversally in the energy surface, along an orbit
which is homoclinic to the torus.
\end{itshape}

The first point is that the tori remain hyperbolic,
and that the 
stable and unstable manifolds are  deformed, but not destroyed
by the additional term.
This results from the observation
that the manifold $M$ formed by the union of the invariant
tori is normally hyperbolic in its energy surface.
Note that this step does not require exponential smallness 
of $\mu$.

It is then a very general result that the stable and unstable manifolds
have non-empty intersection. It is a global property,
which can be established by variational methods, and which still
does not rest on exponential smallness of  $\mu$.

The key point, where exponential smallness is required, 
is transversality.
Since transversality is a generic phenomenon, one may think that
this step is not so crucial.
And indeed, it is very likely that the statement
remains true for most values of $\mu\in ]0,\epsilon]$ (and not only for 
$\mu\leq \mu_0(\epsilon)$).
However, there are two important issues here. First, 
transversality is hard to establish on explicit examples.
Second, it is useful for many further discussions
to obtain some quantitative estimates.

Indeed, we can associate to the intersection between 
the stable and unstable manifolds
a quantity, the splitting, which in a sense measures
transversality. Discussions on such a definition are 
available in \cite{LMS}.
Using methods of Poincar\'e and  Melnikov,
Arnold showed that this splitting can be estimated, 
for sufficiently small $\epsilon$,
by 
\begin{equation}\label{split}
\alpha\geq \mu e^{-C/\sqrt{\epsilon}}+O(\mu^2).
\end{equation}
This implies non-nullity of the splitting, hence 
transversality, for small $\mu$.

\subsection{Transition chain}
We have established the existence, when $\mu>0$ is small enough,
of a family $T_{\omega}$ of hyperbolic invariant tori
such that the stable manifold  $W^+_{\omega}$
and the unstable manifold  $W^-_{\omega}$ intersect transversally
along a homoclinic orbit (but not along $T_{\omega}$!)
for each  $\omega$.

A stability argument shows that the stable manifold 
 $W^+_{\omega}$ of the torus  $T_{\omega}$
intersects transversally the stable manifold
 $W^-_{\omega_0}$ of the torus  $T_{\omega_0}$
when  $\omega$ is close enough to $\omega_0$.
How close directly depends on the size of the splitting.
We obtain heteroclinic orbits between tori close to each other.

Given two values $\omega$ and $\omega'$, we can
find a sequence $\omega_i,1\leq i\leq N$
such that $\omega_0=\omega$, $\omega _N=\omega'$,
and $W^-_i$ intersects transversally $W^+_{i+1}$
for all $i$.
The associated family  $T_{\omega_i}$ of tori is called
a transition chain.

The left step consists in proving that some
orbits shadow the transition chain.
Arnold solved this step by a very simple topological
argument which, however, does not provide any
estimate on the time $T$.
He proves the existence of an orbit joining
any neighborhood of  $T_{\omega}$ to any neighborhood of 
 $T_{\omega'}$.
This ends the proof of the main theorem, since we can chose 
$\omega$ and $\omega'$ such that 
 $\omega<A<B<\omega'$.

The dynamics associated to hyperbolic tori
and transition chains have later been studied more carefully.
It particular, a $\lambda$-lemma can be proved in this
context, which allows to conclude that, in a transition chain,
the unstable manifold  $W^-_0$ of the first torus
intersects transversally the stable manifold of the last torus
$W^+_N$.
These detailed study also allow to relate the speed of diffusion
to the splitting of the invariant manifolds.

\subsection{Diffusion speed}
It is interesting to estimate the speed of evolution
of the variable $p_1$, or in other words the time $T$ in the
statement.
It follows from Nekhoroshev theory that 
this time $T$ has to be exponentially large 
as a function of $\epsilon$.
In fact, it is possible to prove, either by recent developments
on the ideas of Arnold exposed above, or more easily by variational 
methods, \cite{UB},
that
$$
T\leq \frac{e^{C/\sqrt{\epsilon}}}{-\mu \log \mu}
$$
for $\mu\leq \mu_0(\epsilon)$.
This time is of course highly related to the estimate (\ref{split})
of the splitting.
In addition, Ugo Bessi proved that one can take 
$\mu_0(\epsilon)=e^{-C/\sqrt{\epsilon}}$.
Plugging this value of $\mu$ in the estimate of $T$,
we get the estimate 
$T\leq e^{C/\sqrt{\epsilon}}$
as a function of the only parameter $\epsilon$.
Considering 
the fact that the orbit we have described go close to 
double resonances, this is the best estimate one may hope for
in view of the improved Nekhoroshev stability estimates at resonances.

The idea is now well spread that the time of diffusion
is exponentially large. However, we point out that,
if it is indeed  exponentially small
as a function of the parameter $\epsilon$, it is only polynomially
small as a function of the second parameter $\mu$, as was first
understood by P. Lochak and proved in \cite{PB} using the variational
method of U. Bessi.

\subsection{Conclusion}
The theories of instability are developing in several directions.
One of them is to try to understand the limits of stability, and to test
to what extent the stability results obtained so far are optimal.
This aspect has quickly developed recently, for example the
optimal stability exponent $a$ for convex systems is 
almost known.
Another direction is to try to give a description
of unstable orbits in typical systems. This remains a 
widely open question. 

Let us finally mention that the application of the theories
we have presented to concrete systems is very difficult.
One of the reasons is that the  estimates of the threshold
$\epsilon_0$ of validity of Nekhoroshev and KAM theorems
that can (painfully) be obtained by inspection in the proofs
are very bad. Much too bad, for example, to think about
applications to the solar systems with the physical values
of the parameters.

\end{document}